\documentclass[12pt]{article}

\usepackage{mathrsfs}
\usepackage{color}
\usepackage{amssymb, amsthm, amsfonts, amsxtra, amsmath}
\usepackage{mathabx}

\usepackage{latexsym}
\usepackage[all]{xy}
\usepackage{graphics}
\usepackage{latexsym}
\usepackage{makeidx}
\usepackage{rotating}

\theoremstyle{change}

\newtheorem{Theorem}{Theorem}[section]
\newtheorem{Def}[Theorem]{Definition}
\newtheorem{Lem}[Theorem]{Lemma}
\newtheorem{Prop}[Theorem]{Proposition}

\date{}

\begin{document}

\hyphenation{Wo-ro-no-wicz}

\title{Comonoidal W$^*$-Morita equivalence for von Neumann bialgebras}
\author{Kenny De Commer\footnote{Supported in part by the ERC Advanced Grant 227458
OACFT ``Operator Algebras and Conformal Field Theory" }\\ \small Dipartimento di Matematica,  Universit\`{a} degli Studi di Roma Tor Vergata\\
\small Via della Ricerca Scientifica 1, 00133 Roma, Italy\\ \\ \small e-mail: decommer@mat.uniroma2.it}
\maketitle

\newcommand{\acnabla}{\nabla\!\!\!{^\shortmid}}
\newcommand{\undersetmin}[2]{{#1}\underset{\textrm{min}}{\otimes}{#2}}
\newcommand{\otimesud}[2]{\overset{#2}{\underset{#1}{\otimes}}}
\newcommand{\qbin}[2]{\left[ \begin{array}{c} #1 \\ #2 \end{array}\right]_{q^2}}

\abstract{\noindent A theory of Galois co-objects for von Neumann bialgebras is introduced. This concept is closely related to the notion of \emph{comonoidal W$^*$-Morita equivalence} between von Neumann bialgebras, which is a Morita equivalence taking the comultiplication structure into account. We show that the property of `being a von Neumann algebraic quantum group' (i.e.~ `having invariant weights') is preserved under this equivalence relation. We also introduce the notion of a \emph{projective corepresentation} for a von Neumann bialgebra, and show how it leads to a construction method for Galois co-objects and comonoidal W$^*$-Morita equivalences.}

\section*{Introduction}

\noindent In the literature, there are several equivalent ways of introducing the concept of a \emph{W$^*$-Morita equivalence} between von Neumann algebras, for example by means of the categorical formalism (\cite{Rie1}), Connes' correspondences (\cite{Brou1}), Paschke's Hilbert W$^*$-modules (\cite{Pas1}) (called \emph{rigged modules} in \cite{Rie1}), or linking von Neumann algebras (\cite{Bro1}). The latter two approaches will be the ones we favor in this paper. \\

\noindent Let us state the definition of W$^*$-Morita equivalence in terms of linking von Neumann algebras.

\begin{Def} \label{Deflink} (\cite{Bro1},\cite{Rie2}) Let $P$ and $M$ be two von Neumann algebras. A \emph{linking von Neumann algebra} between $P$ and $M$ consists of a von Neumann algebra $Q$ together with a self-adjoint projection $e\in Q$ and $^*$-isomorphisms $P\rightarrow eQe$ and $M\rightarrow (1-e)Q(1-e)$, such that both $e$ and $(1-e)$ are full projections (i.e. have central support equal to 1).\\

\noindent Two von Neumann algebras $P$ and $M$ are called \emph{W$^*$-Morita equivalent} if there exists a linking von Neumann algebra between them.\end{Def}

\noindent In this paper, we will introduce a notion of \emph{comonoidal W$^*$-Morita equivalence} between von Neumann bialgebras. Let us first recall the definition of the latter structure.

\begin{Def}\label{Defvnaqg} A \emph{von Neumann bialgebra} $(M,\Delta_M)$ consists of a von Neumann algebra $M$ and a faithful normal unital $^*$-homomorphism $\Delta_M: M\rightarrow M\bar{\otimes} M$ satisfying the coassociativity condition \[(\Delta_M\otimes \iota)\Delta_M = (\iota\otimes \Delta_M)\Delta_M.\]
\end{Def}

\noindent \emph{Remark:} In the literature, von Neumann bialgebras appear under the name `Hopf-von Neumann algebras'. We prefer to use the above terminology since it is in better correspondence with the purely algebraic nomenclature.  \\

\noindent Our proposal for a notion of comonoidal W$^*$-Morita equivalence between von Neumann bialgebras is the following.

\begin{Def} Let $(P,\Delta_P)$ and $(M,\Delta_M)$ be two von Neumann bialgebras. A \emph{linking weak von Neumann bialgebra} between $(P,\Delta_P)$ and $(M,\Delta_M)$ consists of a linking von Neumann algebra $(Q,e)$ between $P$ and $M$, together with a (non-unital) coassociative normal $^*$-homomorphism $\Delta_Q:Q\rightarrow Q\bar{\otimes} Q$ satisfying \[\Delta_Q(e) = e\otimes e,\qquad \Delta_Q(1-e)= (1-e)\otimes (1-e)\] and, with $Q_{11}=eQe$ and $Q_{22}=(1-e)Q(1-e)$, \[(Q_{11},(\Delta_Q)_{\mid Q_{11}})\cong (P,\Delta_P),\] \[(Q_{22},(\Delta_Q)_{\mid Q_{22}})\cong (M,\Delta_M),\] by the isomorphisms appearing in the definition of a linking von Neumann algebra.\\

\noindent Two von Neumann bialgebras $(P,\Delta_P)$ and $(M,\Delta_M)$ are called \emph{comonoidally W$^*$-Morita equivalent} if there exists a linking weak von Neumann bialgebra between them.
\end{Def}

\noindent We will give some more information on the terminology we use at the beginning of the second section.\\

\noindent In contexts where linking structures appear, one often has a `unilateral version' accompanying it. This one-sided version should then arise as the \emph{corner} of some linking structure. For von Neumann algebras, we will call this structure a Morita Hilbert W$^*$-module (there seems to be no special nomenclature in the literature).

\begin{Def}\label{DefHilb}(\cite{Pas1}) Let $M$ be a von Neumann algebra. A \emph{self-dual (right) Hilbert W$^*$-module} for $M$ consists of a right $M$-module $N$, together with a (non-degenerate) $M$-valued Hermitian inner product $\langle \,\cdot\,,\,\cdot\, \rangle_M$, such that for any bounded $M$-module map $T$ from $N$ to $M$, there exists $x\in N$ for which $T(y)=\langle x,y\rangle_M$ for all $y\in N$.\\

\noindent When the self-dual Hilbert W$^*$-module is \emph{full} (or \emph{saturated}), in the sense that the linear span of all $\langle x,y\rangle_M$, with $x,y\in N$, is $\sigma$-weakly dense in $M$, we call $N$ a \emph{(right) Morita Hilbert W$^*$-module for $M$} (or a Morita Hilbert $M$-module).\\

\end{Def}

\noindent The following definition will then correspond to the unilateral version of a linking weak von Neumann bialgebra.

\begin{Def}\label{Defcoob} A right Galois co-object for a von Neumann bialgebra $(M,\Delta_M)$ consists of a Morita Hilbert W$^*$-module $N$ for $M$, together with a coassociative normal and faithful linear map $\Delta_N:N\rightarrow N\bar{\otimes}N$ for which the following conditions are satisfied:
\begin{enumerate} 
\item For $x\in N$ and $m\in M$, we have $\Delta_N(xm) = \Delta_N(x)\Delta_M(m)$,
\item For $x,y\in N$, we have $\Delta_M(\langle x,y\rangle_M) = \langle \Delta_N(x),\Delta_N(y)\rangle_{M\bar{\otimes} M}$,
\item The linear span of $\{\Delta_N(x)(m_1\otimes m_2)\mid x\in N, m_1,m_2\in M\}$ is $\sigma$-weakly dense in $N\bar{\otimes} N$.
\end{enumerate}\end{Def}

\noindent So the first two conditions give compatibility relations between $\Delta_N,\Delta_M$ and $\langle\,\cdot\,,\,\cdot\,\rangle_M$, while the
final one is a non-degeneracy condition. \\

\noindent Given a notion of Morita equivalence, it is important to consider what properties are invariant under it. The main theorem of this paper will consist of establishing one such an invariant. Let us first introduce the relevant terminology.
 
\begin{Def} (\cite{Kus2},\cite{VDae1}) Let $(M,\Delta_M)$ be a von Neumann bialgebra. We call $(M,\Delta_M)$ a \emph{von Neumann algebraic quantum group} if there exist nsf (normal semi-finite faithful) weights $\varphi_M$ and $\psi_M$ on $M$ such that for all normal states $\omega$ on $M$ and all $x\in M^+$ we have \[\varphi_M((\omega\otimes \iota)\Delta_M(x)) = \varphi_M(x) \qquad \textrm{(left invariance)},\] \[\psi_M((\iota\otimes \omega)\Delta_M(x)) = \psi_M(x) \qquad \textrm{(right invariance)}.\]\end{Def}

\noindent Note that `being a von Neumann algebraic quantum group' is introduced as a property of a von Neumann bialgebra. However, since the weights $\varphi_M$ and $\psi_M$ above turn out to be unique up to scaling with a positive constant, it is customary to consider them as part of the given data.\\
 
\noindent Such von Neumann algebraic quantum groups turn out to have a very rich structure, and seem to form the right framework in which to study the theory of \emph{locally compact quantum groups}. See for example \cite{Kus3}, \cite{Kus2}, \cite{Vae1}, \cite{Kus4}, \cite{Vae2} and \cite{Vae3} for some generalizations to this setting of a large part of the theory of locally compact groups.\\

\noindent The following is the main result of the present paper which we alluded to.

\begin{Theorem}\label{TheovNa} If $(P,\Delta_P)$ and $(M,\Delta_M)$ are comonoidally W$^*$-Morita equivalent von Neumann bialgebras, then $(M,\Delta_M)$ is a von Neumann algebraic quantum group iff $(P,\Delta_P)$ is a von Neumann algebraic quantum group.\end{Theorem}
 
\noindent The proof of this theorem will consist in making the connection with the theory of \cite{DeC1}. Indeed, there a notion of \emph{Galois objects} was introduced. Although one can in fact obtain a complete duality theory between Galois objects (for a von Neumann algebraic quantum group) and Galois co-objects (for the \emph{dual} von Neumann algebraic quantum group), we have refrained from carrying out this discussion in full here, as the details are somewhat technical (in essence, the details of the duality construction can be found in \cite{DeC2}, but one first needs to prove Theorem \ref{TheovNa} of the present paper to be able to use those results).\\
 
\noindent An essential ingredient which allows us to use the theory of \cite{DeC1} will be the notion of a \emph{projective corepresentation} of a von Neumann bialgebra. This notion was also introduced in \cite{DeC1}, but only for von Neumann algebraic quantum groups.
 
\begin{Def}\label{DefProj} Let $(M,\Delta_M)$ be a von Neumann bialgebra. A \emph{(unitary) projective (left) corepresentation} of $(M,\Delta_M)$ on a Hilbert space $\mathscr{H}$ is a left coaction of $(M,\Delta_M)$ on $B(\mathscr{H})$, i.e.~ a normal faithful unital $^*$-homomorphism \[\alpha: B(\mathscr{H})\rightarrow M\bar{\otimes} B(\mathscr{H})\] satisfying the coaction property \[(\iota\otimes \alpha)\alpha = (\Delta_M\otimes \iota)\alpha.\] \end{Def}

\noindent In the third section, we will show that from any projective corepresentation for a von Neumann bialgebra, one can construct from it a Galois co-object for this von Neumann bialgebra. This will generalize the construction of a 2-cocycle function from a projective representation of a (locally compact) group. \\

\noindent As linking von Neumann bialgebras between von Neumann algebraic quantum groups turn out to have a lot of extra structure, such as an associated C$^*$-algebraic description (see again \cite{DeC2}), we prefer to use the following terminology in this case.

\begin{Def}\label{DefvNalqg} Let $(M,\Delta_M)$ and $(P,\Delta_P)$ be von Neumann algebraic quantum groups. Then a linking weak von Neumann bialgebra $(Q,e,\Delta_Q)$ between $(P,\Delta_P)$ and $(M,\Delta_M)$ will be called \emph{a von Neumann algebraic linking quantum groupoid}.\end{Def}

\noindent Indeed, it is intuitively very helpful to see such a von Neumann algebraic linking quantum groupoid between $(P,\Delta_P)$ and $(M,\Delta_M)$ as a kind of $\mathscr{L}^{\infty}$-space on a `quantum groupoid' having a classical object space consisting of two objects, for which the $(M,\Delta_M)$ and $(P,\Delta_P)$ then play the role of `group von Neumann algebras of the isotropy groups', and for which the off-diagonal corners $eQ(1-e)$ and $(1-e)Qe$ play the role of a certain topological linearization of `the space of arrows between the two objects'. See the first section of \cite{DeC5} for some more information (and, for a similar interpretation in a more algebraic setting, see \cite{Bic1}). We note that such von Neumann algebraic linking quantum groupoids then fit into the theory of `measured quantum groupoids' as introduced in \cite{Les1}.\\

\noindent The concrete structure of this paper is as follows.\\

\noindent In \emph{the first section}, we will give some more preliminary information on the notions of linking von Neumann algebras and Hilbert W$^*$-modules.\\

\noindent In \emph{the second section}, we will show how any linking weak von Neumann bialgebra gives rise to a Galois co-object, and, conversely, how any Galois co-object can be completed to a linking weak von Neumann bialgebra. We also show that comonoidal W$^*$-Morita equivalence is indeed an equivalence relation. We end by introducing, in the setting of Galois co-objects for von Neumann algebraic quantum groups, an analogue of the right regular corepresentation.\\

\noindent In \emph{the third section}, we prove the main result concerning projective corepresentations which we mentioned above, and use it to give a proof of Theorem \ref{TheovNa}.\\

\noindent In \emph{the} short \emph{fourth section} we will consider again the special situation of unitary 2-cocycles for a von Neumann bialgebra, which was also treated partly in \cite{DeC1}. Such 2-cocycles correspond precisely to those linking weak von Neumann bialgebras whose underlying linking von Neumann algebra is trivial. We note that, in the operator theoretic framework, these 2-cocycles were introduced in \cite{Eno2}.

\section{W$^*$-Morita equivalence}

\noindent The results in this section are well-known, and most of them are essentially rephrasings of the results in \cite{Pas1}, \cite{Rie1} and \cite{Tak1} (section IX.3). We therefore refrain from giving detailed proofs, but will mostly simply point to the relevant statements in these references.

\subsection{Morita Hilbert W$^*$-modules}

\noindent In Definition \ref{Deflink}, we already recalled what we mean by a linking von Neumann algebra $(Q,e)$ between two von Neumann algebras $P$ and $M$. Let us give some more information on the notation we will use for this concept. First of all, we will always simply identify $P$ and $M$ with their parts inside a linking von Neumann algebra, thus neglecting the identifying maps. We will also write $Q_{ij} = e_{ii}Qe_{jj}$ with $e_{11} = e$ and $e_{22}=1-e$, and \[ Q = \left(\begin{array}{ll} Q_{11} & Q_{12} \\ Q_{21} & Q_{22}\end{array}\right).\] This matrix algebra notation is very convenient in practice. Note that this decomposition makes sense for any projection $e\in Q$, but the special (and characterizing) property of linking von Neumann algebras is that $Q_{12}\cdot Q_{21}$ is $\sigma$-weakly dense in $Q_{11}$ (by definition of fullness for $1-e$), while $Q_{21}\cdot Q_{12}$ is $\sigma$-weakly dense in $Q_{22}$ (by definition of fullness for $e$).\\

\noindent We will further talk simply of `a linking von Neumann algebra' (without specifying what the corners are) or of `a linking von Neumann algebra for the von Neumann algebra $M$' (without specifying the von Neumann algebra in the upper left corner; admittedly, this puts the lower left corner in a privileged position terminology-wise). In fact, this terminology dictates the strongness of the isomorphism one is interested in (keeping none, one or both of the diagonal entries pointwise fixed). The same remark then applies to more general morphisms: if for example $M_1$ and $M_2$ are two von Neumann algebras, $(Q_1,e)$ and $(Q_2,f)$ linking von Neumann algebras for resp.~ $M_1$ and $M_2$, and $\phi:M_1\rightarrow M_2$ a normal unital $^*$-homomorphism, then a $\phi$-compatible unital morphism between $(Q_1,e)$ and $(Q_2,f)$ is a normal unital $^*$-homomorphism $\Phi:Q_1\rightarrow Q_2$ sending $e$ to $f$, whose restriction to a map $Q_{1,22}=M_1\rightarrow Q_{2,22}=M_2$ coincides with $\phi$.\\

\noindent We also defined already the notion of a Morita Hilbert W$^*$-module (Definition \ref{DefHilb}). We introduce the following terminology concerning maps between Morita Hilbert W$^*$-modules.

\begin{Def} When $M_1,M_2$ are two von Neumann algebras, $\phi:M_1\rightarrow M_2$ a unital normal $^*$-homomorphism, and $N_1$ and $N_2$ Morita Hilbert W$^*$-modules for resp.~ $M_1$ and $M_2$, we call a linear map $\Phi:N_1\rightarrow N_2$ a \emph{$\phi$-compatible morphism} when $\Phi(xm)=\Phi(x)\phi(m)$ and $\langle \Phi(x),\Phi(y)\rangle_{M_2} = \phi(\langle x,y\rangle_{M_1})$ for all $x,y\in N_1$ and $m\in M_1$.\\

\noindent When $M$ is a von Neumann algebra, and $N_1$ and $N_2$ two Morita Hilbert $M$-modules, then we call $N_1$ and $N_2$ isomorphic if there exists a bijective $\iota_{M}$-compatible morphism $N_1\rightarrow N_2$, where $\iota_M:M\rightarrow M$ is the identity map.
\end{Def}

\noindent Let us recall from \cite{Pas1}, Proposition 3.10, that if $M$ is a von Neumann algebra, and $N$ a right (Morita) Hilbert $M$-module, then any bounded right $M$-module map $N\rightarrow N$ is adjointable, and the $^*$-algebra of all such maps is a von Neumann algebra. We then introduce the following concept (see \cite{Rie1}).

\begin{Def} (\cite{Rie1}) If $M$ and $P$ are von Neumann algebras, a \emph{$P$-$M$-equivalence bimodule} is a $P$-$M$-bimodule $N$ which is at the same time a right Morita Hilbert $M$-module and left Morita Hilbert $P$-module, and such that \[x\cdot\langle y,z\rangle_M = \langle x,y\rangle_P\cdot z,\qquad \textrm{for all }x,y,z\in N.\]\end{Def}

\noindent The following lemma makes the connection between Morita Hilbert W$^*$-modules and linking von Neumann algebras concrete.

\begin{Lem}\label{Lemid} \begin{enumerate}\item Let $(Q,e)$ be a linking von Neumann algebra between the von Neumann algebras $P$ and $M$. Then $Q_{12}$, together with the $M$-valued inner product \[\langle x,y\rangle_M = x^*y,\qquad x,y\in Q_{12}\] and the $P$-valued inner product \[\langle x,y\rangle_P = xy^*, \qquad x,y\in Q_{12}\] is a $P$-$M$-equivalence bimodule.

\item If $N$ is a right Morita Hilbert $M$-module, there exists a linking von Neumann algebra $(Q,e)$ and an isomorphism $\pi$ of right Hilbert W$^*$-modules from $N$ to $Q_{12}$. Moreover, $(Q,e)$ is then unique up to isomorphism of linking von Neumann algebras for $M$.
\end{enumerate}
\end{Lem}

\begin{proof} The first part of this Lemma can be deduced from Theorem 6.5 of \cite{Rie1}, choosing a concrete representation of $Q$. As for the second part, we can construct the $(Q,e)$ associated to $N$ in a natural way as the von Neumann algebra of right $M$-module maps on the direct sum right Hilbert $W^*$-module $\left(\begin{array}{l} N \\ M\end{array}\right)$ over $M$, together with the projection $e$ onto $N$. The fact that this is then a linking von Neumann algebra follows from the proof of Corollary 7.10 in \cite{Rie1}, which shows that $N$ is a $P$-$M$-equivalence bimodule. Finally, the uniqueness statement follows from Proposition 7.6 of \cite{Rie1}, which shows that in any linking von Neumann algebra $(Q,e)$, the von Neumann algebra $Q_{11}$ can be identified with the set of bounded right $Q_{22}$-module maps on $Q_{12}$. This then easily allows one to identify this linking von Neumann algebra with the canonical one we constructed above.
\end{proof}

\noindent In the following, we will always regard a Morita Hilbert W$^*$-module as the upper right corner of its associated linking von Neumann algebra. This allows us to introduce a lot of operations for Morita Hilbert W$^*$-modules in a straightforward way. For example, if $N$ is a Morita Hilbert W$^*$-module, and $(Q,e)$ the associated linking von Neumann algebra, then the \emph{predual} of $N$, whose existence was proven in \cite{Pas1}, may be identified with the space of normal functionals on $Q$ which vanish on all $Q_{ij}$ except $Q_{12}$. The $\sigma$-weak topology of $N$ as the dual of its predual then coincides with the restriction of the $\sigma$-weak topology on $N\subseteq Q$. This allows us to talk about normal maps between Morita Hilbert W$^*$-modules without any ambiguity.\\

\noindent The following Lemma shows how to complete maps which are only defined on a subspace of a Morita Hilbert W$^*$-module.

\begin{Lem}\label{LemuniMod} Let $M_1,M_2$ be von Neumann algebras, equipped with a unital normal $^*$-homomorphism $\phi:M_1\rightarrow M_2$. Let $N_1$ and $N_2$ be right Morita Hilbert W$^*$-modules over resp.~ $M_1$ and $M_2$. Suppose that $\mathscr{N}_1$ is a $\sigma$-weakly dense $M$-submodule of $N_1$, and suppose that there exists a linear map $\pi: \mathscr{N}_1\rightarrow N_2$ such that $\pi(xm)=\pi(x)\phi(m)$ and $\langle \pi(x),\pi(y)\rangle_{M_2} = \phi(\langle x,y\rangle_{M_1})$ for all $x,y\in \mathscr{N}_1$ and $m\in M$. Then $\pi$ has a unique extension to a normal $\phi$-compatible morphism $\Psi:N_1\rightarrow N_2$. If $\phi$ is faithful, then $\Psi$ will be faithful. If $\phi$ is bijective, and $\pi$ has $\sigma$-dense image, then $\Psi$ is bijective.
\end{Lem}

\begin{proof} As $\mathscr{N}_1$ is a linear space, it is also $\sigma$-strongly dense in $N_1$. Further, from the $\phi$-compatibility condition on $\pi$, we easily get that if a net $x_\alpha\in \mathscr{N}_1$ converges $\sigma$-strongly to 0, then also $\pi(x_{\alpha})\rightarrow 0$ in the $\sigma$-strong topology. From these two observations, it follows that $\pi$ can be uniquely extended to a normal map $\Psi:N_1\rightarrow N_2$, which is then of course still $M$-linear and $\phi$-compatible.\\

\noindent If $\phi$ is faithful, then $\Psi(x)=0$ for $x\in N_1$ would imply $\phi(\langle y,x\rangle_{M_1})=0$ for all $y\in \mathscr{N}_1$, hence $x=0$; thus also $\Psi$ is faithful. If further $\phi$ is bijective and $\pi$ has $\sigma$-dense image, then, as the range $\Psi(N_1)$ is $\sigma$-weakly closed, it must equal $N_2$, and hence $\Psi$ is bijective.

\end{proof}

\noindent The next Lemma provides a further weakening of the conditions in the previous Lemma.

\begin{Lem}\label{LemuniMod2} Let $M_1,M_2$ be von Neumann algebras, $\phi:M_1\rightarrow M_2$ a unital normal $^*$-homomorphism. Let $N_1, N_2$ be right Morita Hilbert W$^*$-modules for resp.~ $M_1$ and $M_2$. Let $I$ be an index set, and suppose $x_i \in N_1$ and $y_i\in N_2$ are elements such that $\phi(\langle x_i,x_j\rangle_{M_1}) = \langle y_i,y_j\rangle_{M_2}$ for all $i,j\in I$, and suppose that the $M_1$-linear span of the $x_i$ is $\sigma$-weakly dense in $N_1$. Then there exists a unique $\phi$-compatible morphism $\pi:N_1\rightarrow N_2$ of Hilbert W$^*$-modules such that $\pi(x_i)=y_i$.\end{Lem}

\begin{proof} Let $\mathscr{N}_1$ be the right $M_1$-module spanned by the $x_i$. Then the map \[\pi:\mathscr{N}_1\rightarrow N_2: \sum_{i=1}^n x_im_i\rightarrow \sum_{i=1}^n y_i\phi(m_i),\qquad m_i\in M\] is a well-defined $\phi$-intertwining map, since, by the compatibility between the $x_i$ and $y_i$, we have \[\langle \sum_{i=1}^n x_im_i,\sum_{i=1}^n x_im_i\rangle_{M_1}=0 \qquad\Rightarrow\qquad \langle \sum_{i=1}^n y_i\phi(m_i),\sum_{i=1}^n y_i\phi(m_i)\rangle_{M_2}=0.\] The Lemma then follows immediately by the previous one.
\end{proof}

\noindent The way in which linking von Neumann algebras most frequently appear is the following (see also Theorem 8.15 and its footnote in \cite{Rie1}). The proof of the Proposition essentially follows by Proposition 1.3 and Proposition 1.1.(2) of \cite{Rie1}.

\begin{Prop}\label{PropLinkk} Let $Z$ be a von Neumann algebra, and let $\mathscr{H}_1$ and $\mathscr{H}_2$ be two Hilbert spaces equipped with \emph{faithful} normal $^*$-representations $\pi_1$ and $\pi_2$ of $Z$. Denote $M = \pi_2(Z)'$ and $P=\pi_1(Z)'$. Then the space $N$ of $\pi_1$-$\pi_2$-intertwiners is a right Morita Hilbert $M$-module, and the commutant $Q$ of the direct sum representation $\pi_1\oplus \pi_2$, together with the projection $e$ on $\mathscr{H}_1$, is a linking von Neumann algebra between $P$ and $M$.
\end{Prop}

\noindent In particular, this shows that the notion of `linking algebra' which was used in \cite{DeC1} coincides with the terminology of the present paper. \\

\noindent Another way to create Morita Hilbert W$^*$-modules is the following. It is essentially a concrete, spatial approach to \emph{ternary W$^*$-algebras}.

\begin{Prop}\label{PropTern} Let $\mathscr{H}$ and $\mathscr{K}$ be two Hilbert spaces, and let $N \subseteq B(\mathscr{H},\mathscr{K})$ be a $\sigma$-weakly closed linear space for which the set $\{xy^*z\mid x,y,z\in N\}$ equals $N$. Then with $M$ denoting the $\sigma$-weak closure of the linear span of $\{x^*y\mid x,y\in N\}$, we have that $M$ is a von Neumann algebra and $N$ a right Morita Hilbert $M$-module for the $M$-valued inner product $\langle x,y\rangle_M = x^*y$.
\end{Prop}

\begin{proof} Denote $O = N^*$, the set of adjoints of elements in $N$. By the condition $\{xy^*z\mid x,y,z\in N\}=N$, we have that $O\cdot N=\{\sum_{i=1}^n x_i^*y_i\mid n\in \mathbb{N}_0,x_i,y_i\in N\}$ and $N\cdot O$ are $^*$-algebras. Hence their respective $\sigma$-weak closures $M$ and $P$ are von Neumann algebras (possibly with different units than $1_{B(\mathscr{H})}$ and $1_{B(\mathscr{K})}$). As $N$ is $\sigma$-weakly closed, $N$ is a $P$-$M$-bimodule, and then it is immediate that $\left(\begin{array}{cc} P & N \\O & M\end{array}\right)$ is a von Neumann algebra. By the way $M$ and $P$ were defined, it is a linking von Neumann algebra between $P$ and $M$. In particular, $N$ is a right Morita Hilbert $M$-module.
\end{proof}

\noindent We also record the following Lemma for further use.

\begin{Lem}\label{Lemext} Let $M_1$ and $M_2$ be von Neumann algebras, and $N_1$ and $N_2$ Morita Hilbert W$^*$-modules for resp.~ $M_1$ and $M_2$. Let $\pi_{22}:M_1\rightarrow M_2$ be a normal unital $^*$-homomorphism, and $\pi_{12}:N_1\rightarrow N_2$ a $\pi_{22}$-compatible normal morphism  Then if $(Q_1,e)$ and $(Q_2,f)$ are the linking von Neumann algebras associated to respectively $N_1$ and $N_2$, there exists a unique $\pi_{22}$-compatible, \emph{not necessarily unital morphism} $\pi:(Q_1,e)\rightarrow (Q_2,f)$ such that $\pi(e)\leq f$, $\pi(1-e) = 1-f$, and such that the restriction to $N_1$ coincides with $\pi_{12}$.\\

\noindent If the right $M_2$-module generated by $\pi_{12}(N_1)$ is $\sigma$-weakly dense in $N_2$, then $\pi(e)=f$, and hence $\pi$ unital.\end{Lem}

\begin{proof} The uniqueness of $\pi$ is immediate. Also the existence of $\pi:(Q_1,e)\rightarrow (Q_2,f)$ as a normal $^*$-homomorphism follows from basic von Neumann algebraic techniques. If $\pi_{12}(N_1)\cdot M_2$ is $\sigma$-weakly dense in $N_2$, then $\pi(e)$ acts as a unit on $N_2$ by left multiplication, and hence equals $f$.
\end{proof}

\subsection{Tensor products and composition}\label{SubsecTen}

\noindent Suppose that $M_1,M_2,P_1,P_2$ are von Neumann algebras, and that $(Q_1,e)$, resp.~ $(Q_2,f)$, is a linking von Neumann algebra between $P_1$ and $M_1$, resp.~ $P_2$ and $M_2$. Then we denote $Q_1* Q_2$ for the corner of $Q_1\bar{\otimes}Q_2$ by the projection $e\otimes f + (1-e)\otimes (1-f)$. The reason for this notation is that this can (easily) be shown to be a special case of a fibred product of von Neumann algebras (i.c.~ fibred over $\mathbb{C}^2$), see \cite{Eno3}, sections 2.3 and 2.4.\\

\noindent It is easy to see that $(Q_1* Q_2,e\otimes f)$ will be a linking von Neumann algebra between $P_1\bar{\otimes} P_2$ and $M_1\bar{\otimes}M_2$. The operation $*$ is an associative operation on linking von Neumann algebras.\\

\noindent If $N_1$, resp.~ $N_2$, is a Morita Hilbert W$^*$-module for a von Neumann algebra $M_1$, resp.~ $M_2$, we can define $N_1\bar{\otimes} N_2$ to be the right hand corner of $Q_1* Q_2$, with $(Q_1,e)$ and $(Q_2,f)$ the linking von Neumann algebras associated with respectively $N_1$ and $N_2$. We then have a natural injection $N_1\odot N_2 \rightarrow N_1\bar{\otimes}N_2$, where $\odot$ denotes the algebraic tensor product, and this allows us also to see $N_1\bar{\otimes}N_2$ as a concrete realization of the `self-dual completion' of the pre-Hilbert W$^*$-module $N_1\odot N_2$ for $M_1\bar{\otimes}M_2$ (see Theorem 3.2 of \cite{Pas1}, and also Proposition 8.5 in \cite{Rie1}). Moreover, $(Q_1* Q_2,e\otimes f)$ will then be a linking von Neumann algebra associated to the Morita Hilbert W$^*$-module $N_1\bar{\otimes}N_2$ over $M_1\bar{\otimes} M_2$. In the same way, we can take the tensor product of the lower left corners of $Q_1$ and $Q_2$, and thus, if we write $Q_i = \left(\begin{array}{ll} P_i & N_i \\ O_i & M_i\end{array}\right)$, we can write \[Q_1* Q_2 = \left(\begin{array}{cc} P_1\bar{\otimes}P_2 & N_1\bar{\otimes}N_2\\O_1\bar{\otimes}O_2 & M_1\bar{\otimes} M_2\end{array}\right).\] Finally, if $N_1$, $N_2$ and $N_3$ are Morita Hilbert W$^*$-modules, and $\phi:N_1\rightarrow N_2$ a normal map, it is clear, by passing again to the enveloping linking von Neumann algebra picture, that one can define a slice map $\phi\otimes \iota:N_1\bar{\otimes}N_3\rightarrow N_2\bar{\otimes}N_3$, uniquely determined by the property that it is normal and satisfies $(\phi\otimes \iota)(x\otimes y)=\phi(x)\otimes y$ for elementary tensors $x\otimes y \in N_1\bar{\otimes} N_3$.\\

\noindent Let us also comment on how Morita Hilbert W$^*$-modules can be \emph{composed}, which will show in particular that W$^*$-Morita equivalence is an equivalence relation. Let $M_1,M_2$ and $M_3$ be von Neumann algebras, and let $N_{12}$ be an $M_1$-$M_2$-equivalence bimodule, and $N_{23}$ an $M_2$-$M_3$-equivalence bimodule. Consider the associated linking von Neumann algebras, which we will denote as \[Q_1 = \left(\begin{array}{cc} M_{1} & N_{12} \\ N_{21} & M_2\end{array}\right),\qquad Q_2 = \left(\begin{array}{cc} M_2 & N_{23} \\ N_{32} & M_3\end{array}\right).\] Then we can consider the direct sum right Hilbert W$^*$-module $\left(\begin{array}{c} N_{12} \\ M_2 \\ N_{32}\end{array}\right)$ for $M_2$. Let $Q$ be the von Neumann algebra of bounded right $M_2$-linear maps on this module (using again Proposition 3.10 of \cite{Pas1}). Then we can decompose $Q$ as \[Q = \left(\begin{array}{ccc} M_1 & N_{12} & N_{13} \\ N_{21} & M_2 & N_{23} \\ N_{31} & N_{32} & M_3\end{array}\right).\] As then $N_{13}\supseteq N_{12}\cdot N_{23}$ and $N_{31}\supseteq N_{32}\cdot N_{21}$, we see that $N_{31}\cdot N_{13}$ contains $N_{32}\cdot (N_{21}\cdot N_{12})\cdot N_{23}$. As $N_{21}\cdot N_{12}$ is $\sigma$-weakly dense in $M_2$, and $N_{32}\cdot N_{23}$ is $\sigma$-weakly dense in $M_3$, we get that $N_{31}\cdot N_{13}$ is $\sigma$-weakly dense in $M_3$. Similarly, $N_{13}\cdot N_{31}$ is $\sigma$-weakly dense in $M_1$. This implies that \[ \left(\begin{array}{cc} M_1 & N_{13} \\ N_{31} & M_3\end{array}\right)\] is a linking von Neumann algebra between $M_1$ and $M_3$, which we call the \emph{composition} of $Q_1$ and $Q_2$. The $N_{13}$-part, considered as a $M_1$-$M_3$-equivalence bimodule, is called the \emph{composition} of the equivalence bimodules $N_{12}$ and $N_{23}$. One could also call the total structure $Q$, together with the units of its diagonal components, a `3$\times$3 linking von Neumann algebra'.

\subsection{Compatibility with weight theory}\label{SubsecWeight}

\noindent Let us now comment on the relation with weight theory for von Neumann algebras (see \cite{Tak1}, and especially Chapter IX, section 3 for a discussion of material closely related to ours).\\

\noindent Let $M$ be a von Neumann algebra, and let $\psi$ be an nsf (i.e., normal semi-finite faithful) weight on $M$. We denote by $\mathscr{N}_{M,\psi}$ the space of elements $x\in M$ for which $\psi(x^*x)<\infty$, by $\mathscr{M}_{\psi}^{+}$ the space of elements $x\in M^+$ for which $\psi(x)<\infty$, and by $\mathscr{M}_{\psi}$ we denote the linear span of $\mathscr{M}_{\psi}^+$, which also coincides with $\mathscr{N}_{M,\psi}^*\cdot \mathscr{N}_{M,\psi}$. Then one can linearly extend $\psi$ to $\mathscr{M}_{\psi}$, and we will use the same notation for this extension.\\

\noindent Now suppose that $N$ is a right Morita Hilbert $M$-module. Then we can also form the space $\mathscr{N}_{N,\psi}$ of elements $x\in N$ for which $\psi(\langle x,x\rangle_M)<\infty$. Clearly, this space is $\sigma$-weakly dense in $N$ as it contains the set $N\cdot \mathscr{N}_{M,\psi}$ (in fact, it \emph{equals} this set by a polar decomposition argument). We can then turn $\mathscr{N}_{N,\psi}$ into a pre-Hilbert space by the scalar product $\langle x,y\rangle = \psi(\langle x,y\rangle_M)$ (we will take the scalar product in our Hilbert spaces conjugate linear in the first variable, as this is the most natural thing to do in this context). We denote by $\mathscr{L}^2(N,\psi)$ its completion, and by $\Gamma_{N,\psi}$ the natural embedding map $\mathscr{N}_{N,\psi}\hookrightarrow \mathscr{L}^2(N,\psi)$. Applying the same construction to $M$ considered as a right Morita $M$-module, we obtain the ordinary GNS-construction associated to $\psi$. The latter however also comes with a normal left representation $\pi_{M,\psi}$ of $M$ on $\mathscr{L}^2(M)$, uniquely determined by the property that $\pi_{M,\psi}(x)\Gamma_{M,\psi}(y) = \Gamma_{M,\psi}(xy)$ for $y\in \mathscr{N}_{\psi}$ and $x\in M$. A similar left representation can then be obtained for $N$, but it will not act on one Hilbert space, but as linear operators between two different Hilbert spaces. Namely, for $x\in N$ and $y\in \mathscr{N}_{M,\psi}$, we have that $\|\Gamma_{N,\psi}(xy)\| \leq \|\langle x,x\rangle_M\|^{1/2}\,\|\Gamma_{M,\psi}(y)\|$, so that one can define $\pi_{N,\psi}(x)$ as the unique bounded linear operator \[\pi_{N,\psi}(x):\mathscr{L}^2(M,\psi)\rightarrow \mathscr{L}^2(N,\psi)\textrm{ such that }\pi_{N,\psi}(x)\Gamma_{\psi}(y)\rightarrow \Gamma_{\psi}(xy)\textrm{ for all }y\in \mathscr{N}_{M,\psi}.\] Then $\pi_{N,\psi}$ will be a normal map of $N$ into $B(\mathscr{L}^2(M,\psi),\mathscr{L}^2(N,\psi))$, and clearly \[\pi_{N,\psi}(xy)=\pi_{N,\psi}(x)\pi_{M,\psi}(y)\qquad \textrm{ for all }x\in N \textrm{ and }y\in M.\] It is also easily computed that \[\pi_{N,\psi}(x)^*\pi_{N,\psi}(y) = \pi_{M,\psi}(\langle x,y\rangle_M) \textrm{ for }x,y\in N.\]

\noindent If then $(Q,e)$ is the linking von Neumann algebra associated with $N$, we can represent it in a faithful, normal and unit-preserving way on $\left(\begin{array}{cc} \mathscr{L}^2(N,\psi)\\\mathscr{L}^2(M,\psi)\end{array}\right)$, again essentially by extending the left multiplication operation on $\left(\begin{array}{cc} \mathscr{N}_{N,\psi}\\\mathscr{N}_{M,\psi}\end{array}\right)$. In particular, we have a unital faithful normal $^*$-representation of $P$ on $\mathscr{L}^2(N)$. The above constructions can further be brought in connection with the theory of GNS-representations for $Q$, and one could also develop a theory of `standard' representations. However, in this paper, we will not need this further structure, so we refrain from making these further elaborations.\\

\noindent As we will only need one nsf weight at any particular moment, we will in the following unburden the notation somewhat by dropping the symbol $\psi$ in the notation for the GNS-construction.\\

\noindent Let us now give some comments on the tensor product theory of weights. If $M_1$ and $M_2$ are von Neumann algebras, and $\psi_i$ an nsf weight on $M_i$, then one can define the tensor product weight $\psi_1\otimes \psi_2$ on $M_1\bar{\otimes}M_2$. In Definition VIII.4.2 of \cite{Tak1}, this is introduced by using the language of (left) Hilbert algebras. Alternatively, $\psi_1\otimes \psi_2$ can also be introduced using operator valued weights: one can consider $(\iota\otimes \psi_2)$ as an (nsf) operator valued weight from $M_1\bar{\otimes}M_2$ to $M_2$, while $(\psi_1\otimes \iota)$ can be considered an (nsf) operator valued weight from $M_1\bar{\otimes}M_2$ to $M_1$. Then $\psi_1\circ (\iota\otimes \psi_2)$ and $\psi_2\circ (\psi_1\otimes\iota)$ are well-defined nsf weights on $M_1\bar{\otimes}M_2$, and they can be shown to be equal to each other (for example, by using that an nsf weight can be written as the pointwise limit of a net of increasing positive functionals). It can then be shown that this agrees with the nsf weight $\psi_1\otimes \psi_2$ as defined in the first way.\\

\noindent Let now $N_1$ and $N_2$ be right Morita Hilbert W$^*$-modules for respective von Neumann algebras $M_1$ and $M_2$, and $\psi_i$ an nsf weight on $M_i$. Then one can identify $\mathscr{L}^2(N_1\bar{\otimes}N_2)$ unitarily with $\mathscr{L}^2(N_1)\otimes \mathscr{L}^2(N_2)$ by the unique unitary which sends $\Gamma_{N_1\bar{\otimes}N_2}(x\otimes y)$ into $\Gamma_{N_1}(x)\otimes \Gamma_{N_2}(y)$ for $x\in \mathscr{N}_{N_1,\psi_1}$ and $y\in \mathscr{N}_{N_2,\psi_2}$. In the following, we will then always use $\mathscr{L}^2(N_1)\otimes \mathscr{L}^2(N_2)$ for the GNS-space of $\psi_1\otimes \psi_2$, but we will then write the associated GNS-map as $\Gamma_{N_1}\otimes \Gamma_{N_2}$. Of course, the associated representation of $N_1\bar{\otimes}N_2$ then becomes the tensor product representation $\pi_{N_1}\otimes\pi_{N_2}$ into $B(\mathscr{L}^2(M_1)\otimes \mathscr{L}^2(M_2), \mathscr{L}^2(N_1)\otimes \mathscr{L}^2(N_2))$.\\

\section{Comonoidal W$^*$-Morita equivalence}

\noindent Suppose that $P$ and $M$ are von Neumann algebras which also have some extra structure. One would then like an appropriate kind of W$^*$-Morita equivalence which takes this structure into account. This leads quite naturally to the notion of \emph{comonoidal} W$^*$-Morita equivalence between von Neumann bialgebras, introduced in Definition \ref{Deflink}. Let us remark that the notion of a linking weak von Neumann bialgebra $(Q,e,\Delta_Q)$ can also be defined more succinctly using the operation $*$ introduced above in subsection \ref{SubsecTen}. Indeed, then it becomes simply a linking von Neumann algebra $(Q,e)$ equipped with a coassociative normal unital morphism $\Delta_Q: (Q,e)\rightarrow (Q,e)*(Q,e)=(Q* Q,e\otimes e)$. We will further use the following simplifying notation: \[\Delta_{ij}: Q_{ij}\rightarrow Q_{ij}\bar{\otimes}Q_{ij}\] denotes the restriction of $\Delta_Q$ to $Q_{ij}$. We also follow the same conventions as for linking von Neumann algebras, and will talk about `a linking weak von Neumann bialgebra' or `a linking weak von Neumann bialgebra for $(M,\Delta_M)$'. \\

\noindent Let us comment now on the terminology we use. The term `weak von Neumann bialgebra' is a straightforward analogue of the notion of a `weak bialgebra', as introduced in \cite{Boh1}. (Although the terminology von Neumann weak bialgebra would then be more accurate, this seems more awkward to use.) The terminology `von Neumann algebraic linking quantum groupoid' (Definition \ref{DefvNalqg}) has already been motivated somewhat in the introduction. Finally, to explain the terminology `comonoidal', let us suppose for the moment that we are in the finite-dimensional setting, and that we do not consider the associated $^*$-structure. Then it is not difficult to show that if $(Q,e,\Delta_Q)$ is a `linking weak bialgebra', we have an equivalence functor $M$-$\textrm{Mod}\rightarrow P$-$\textrm{Mod}$ by taking the balanced tensor product on the left with $\,\!_PN_M$. This equivalence functor is naturally endowed with a weak comonoidal structure $F$. Namely, if $V,W\in M$-$\textrm{Mod}$, we have \[F: N\underset{M}{\otimes} (V\otimes W)\rightarrow (N\underset{M}{\otimes}V)\otimes (N\underset{M}{\otimes} W): x\underset{M}{\otimes}(v\otimes w)\rightarrow (x_{(1)}\underset{M}{\otimes} v)\otimes (x_{(2)}\underset{M}{\otimes} w),\] where we have used the Sweedler notation for $\Delta_N$. In case the corners of $Q$ are \emph{Hopf} algebras, this weak comonoidal structure can be shown to be strong. A similar discussion then holds in the analytic setting: for a general linking weak von Neumann bialgebra, we will get a weakly comonoidal $^*$-equivalence between the monoidal categories $\textrm{Rep}^*$ of normal unital $^*$-representations of the corner von Neumann algebras on Hilbert spaces, and this will be strongly comonoidal if these corners are von Neumann algebraic quantum groups (see again \cite{DeC2} for details). In any case, we have seen that it is the comonoidal structure which appears most naturally, hence we use it to designate the structure.\\

\noindent In the introduction, we also introduced the notion of a Galois co-object (Definition \ref{Defcoob}). Let us remark that one may drop the assumption of faithfulness and normality of the map $\Delta_N$ in that definition, as they are a consequence of the second compatibility condition.\\

\noindent The following Proposition provides the connection between Galois co-objects the linking weak von Neumann bialgebras.

\begin{Prop}\label{PropCooblink} Let $(N,\Delta_N)$ be a right Galois co-object for a von Neumann bialgebra $(M,\Delta_M)$, and let $(Q,e)$ be a linking von Neumann algebra associated to $N$. Then there exists a unique linking weak von Neumann bialgebra structure $\Delta_Q$ on $(Q,e)$ such that the restriction of $\Delta_Q$ to $N$ coincides with $\Delta_N$.\\

\noindent Conversely, if $(Q,e,\Delta_Q)$ is a linking weak von Neumann bialgebra for a von Neumann bialgebra $(M,\Delta_M)$, then the upper right hand corner $(Q_{12},\Delta_{12})$ is a Galois co-object for $(M,\Delta_M)$.

\end{Prop}

\begin{proof} Let $(N,\Delta_N)$ be a right von Neumann algebraic Galois co-object, and let $(Q,e)= \left(\begin{array}{ll} P & N \\ O & M\end{array}\right)$ be the linking von Neumann algebra associated to $N$ as in Lemma \ref{Lemid}. Then we can apply Lemma \ref{Lemext} with respect to $\Delta_N$ and $\Delta_M$ to obtain a faithful normal $^*$-homomorphism $\Delta_Q: Q\rightarrow Q* Q \subseteq Q\otimes Q$ with $\Delta_Q(1-e)=1-e$. By the uniqueness statement in that Lemma, we have that $\Delta_Q$ is coassociative, since $(\Delta_Q\otimes \iota)\Delta_Q$ and $(\iota\otimes \Delta_Q)\Delta_Q$ coincide when restricted to $N$ and $M$. As $\Delta_N(N)(M\bar{\otimes}M)$ is $\sigma$-weakly dense in $N\bar{\otimes}N$ by definition of a Galois co-object, the `non-degeneracy' condition in that Lemma is satisfied, so that $\Delta_Q:Q\rightarrow Q* Q$ is unital. Hence $(Q,e,\Delta_Q)$ is a linking weak von Neumann bialgebra.\\

\noindent Conversely, suppose that $(Q=\left(\begin{array}{ll} P & N \\ O & M\end{array}\right),\Delta_Q)$ is a linking weak von Neumann bialgebra. Then it is clear that $(N,\Delta_N)$ satisfies the first two conditions of a Galois co-object. Suppose that $\Delta_N(N)(M\bar{\otimes}M)$ is \emph{not} $\sigma$-weakly dense in $N\bar{\otimes} N$. Since the former space is a non-trivial right $M\bar{\otimes}M$-module, we can find a non-zero $x\in P\bar{\otimes}P$ such that $x \Delta_N(y)=0$ for all $y\in N$. (Indeed: then the $\sigma$-weak closure of $\Delta_N(N)(O\bar{\otimes} O)$ will be a non-trivial right ideal inside $P\bar{\otimes}P$, hence there exists a non-zero projection $x\in P\bar{\otimes}P$ which annihilates it by left multiplication.) But then $x\Delta_P(yz)=0$ for all $y\in N, z\in O$. Since the space $N\cdot O$ is $\sigma$-weakly dense in $P$, also $x\Delta_P(w)=0$ for all $w\in P$. Since $\Delta_P(1_P)=1_P\otimes 1_P$, we find that $x=0$, a contradiction. Hence $\Delta_N(N)(M\bar{\otimes}M)$ is $\sigma$-weakly dense in $N\bar{\otimes} N$.
\end{proof}

\noindent \emph{Remark:} If $(M,\Delta_M)$ is a von Neumann algebraic quantum group, we know that $\Delta_M(M)(1\otimes M)$ is $\sigma$-weakly dense in $M\bar{\otimes}M$ (this follows from Corollary 6.11 of \cite{Kus1}, applied to the associated reduced C$^*$-algebraic quantum group)). Hence in this case, we may relax the density condition for a von Neumann algebraic Galois co-object to `$\Delta_N(N)(1\otimes M)$ being $\sigma$-weakly dense in $N\bar{\otimes} N$'. This is more in line with the way Galois co-objects are defined in the setting of Hopf algebras (see \cite{Schn1}, section 4, although the terminology of Galois co-object is not used there).\\

\noindent The following Proposition is mandatory to prove if we want to use the terminology introduced.\\

\begin{Prop} Comonoidal W$^*$-Morita equivalence induces an equivalence relation between von Neumann bialgebras.
\end{Prop}

\begin{proof} It is clear that if $(M,\Delta_M)$ is a von Neumann bialgebra, then it is comonoidally W$^*$-Morita equivalent with itself by the linking weak von Neumann bialgebra $(Q,\Delta_Q)$ which has $Q= M\otimes M_2(\mathbb{C})$, and with $\Delta_{ij}=\Delta_M$ on $Q_{ij}=M$. Further, if $(P,\Delta_P)$ and $(M,\Delta_M)$ are comonoidally W$^*$-Morita equivalent by a linking weak von Neumann bialgebra $(Q,e,\Delta_Q)$, then also $(M,\Delta_M)$ and $(P,\Delta_P)$ are, by the linking weak von Neumann bialgebra $(Q,1-e,\Delta_Q)$.\\

\noindent Now let $(Q_1,e,\Delta_{Q_1})$ and $(Q_2,f,\Delta_{Q_2})$ be two linking weak von Neumann bialgebras. As explained in the second part of Subsection \ref{SubsecTen}, we can combine $(Q_1,e)$ and $(Q_2,f)$ into a global 3$\times$3-linking von Neumann algebra \[ Q = \left(\begin{array}{lll} Q_{11} & Q_{12} & Q_{13} \\ Q_{21} & Q_{22} & Q_{23} \\ Q_{31} & Q_{32} & Q_{33}\end{array}\right),\] with $(Q_1,e)$ isomorphic to the upper left hand block, and $(Q_2,f)$ isomorphic to the lower right hand block. We then have an obvious extension of $*$ to such 3$\times$3-linking von Neumann algebras (which is then a fibred product over $\mathbb{C}^3$), and we can write \[ Q* Q = \left(\begin{array}{lll} Q_{11}\bar{\otimes}Q_{11} & Q_{12}\bar{\otimes} Q_{12} & Q_{13}\bar{\otimes}Q_{13} \\ Q_{21}\bar{\otimes}Q_{21} & Q_{22}\bar{\otimes}Q_{22} & Q_{23}\bar{\otimes} Q_{23} \\ Q_{31}\bar{\otimes}Q_{31} & Q_{32}\bar{\otimes}Q_{32} & Q_{33}\bar{\otimes}Q_{33}\end{array}\right).\]

\noindent Transporting the comultiplication structures from $(Q_1,e,\Delta_{Q_1})$ and $(Q_2,e,\Delta_{Q_2})$, we then have maps $\Delta_{ij}:Q_{ij}\rightarrow Q_{ij}\bar{\otimes} Q_{ij}$ for $|i-j|\leq 1$. Now denote by $\widetilde{Q}_{13}$ the set $Q_{12}\cdot Q_{23}$, which will then be a $\sigma$-weakly dense subset of $Q_{13}$ (since, if not, it would have, being a right $Q_{33}$-module, a non-zero left annihilator in $Q_{11}$, which is clearly impossible as $Q_{12}\cdot Q_{23}\cdot Q_{32}\cdot Q_{21}$ is $\sigma$-weakly dense in $Q_{11}$). Applying Lemma \ref{LemuniMod2} to the elements $xy$ and $\Delta_{12}(x)\Delta_{23}(y)$ for $x\in Q_{12},y\in Q_{23}$, we see that we can find a normal faithful linear map $\Delta_{13}:Q_{13}\rightarrow Q_{13}\bar{\otimes}Q_{13}$, which will then be $\Delta_{11}$-$\Delta_{33}$-compatible and coassociative. Defining \[\Delta_{31}:Q_{31}\rightarrow Q_{31}\bar{\otimes}Q_{31}: x\rightarrow (\Delta_{13}(x^*))^*,\] we get that \[(\left(\begin{array}{ll} Q_{11} & Q_{13} \\ Q_{31} & Q_{33}\end{array}\right),\left(\begin{array}{ll} \Delta_{11} & \Delta_{13} \\ \Delta_{31} & \Delta_{33}\end{array}\right))\] is a linking weak von Neumann bialgebra between $(Q_{11},\Delta_{11})$ and $(Q_{33},\Delta_{33})$. From this, it follows immediately that comonoidal W$^*$-Morita equivalence is a transitive relation, which finishes the proof.
\end{proof}

\noindent We now construct, in the setting of Galois co-objects for von Neumann algebraic quantum groups, an analogue of the right regular corepresentation for a von Neumann algebraic group.\\

\begin{Prop} Let $(M,\Delta_M)$ be a von Neumann algebraic quantum group with a right invariant nsf weight $\psi_M$. Let $(N,\Delta_N)$ be a right Galois co-object for $(M,\Delta_M)$. Then for all $x\in \mathscr{N}_{N,\psi_M}$ and $y\in \mathscr{N}_{M,\psi_M}$, the element $\Delta_N(x)(1\otimes y)$ lies in $\mathscr{N}_{N\bar{\otimes}N,\psi_M\otimes\psi_M}$, and there exists a unitary element $\widetilde{V} \in B(\mathscr{L}^2(N))\bar{\otimes} N$ such that \[\widetilde{V}\,\Gamma_N(x)\otimes \Gamma_M(y) = (\Gamma_N\otimes \Gamma_N)(\Delta_N(x)(1\otimes y)).\]

\noindent Furthermore, if $x\in \mathscr{N}_{N,\psi_M}$ and $\omega\in N_*$, then $(\iota\otimes \omega)(\Delta_N(x))\in \mathscr{N}_{N,\psi_M}$, and \[(\iota\otimes \omega)(\widetilde{V})\Gamma_N(x) = \Gamma_N((\iota\otimes \omega)\Delta_N(x)).\]
\end{Prop}

\begin{proof} The proof that $\widetilde{V}$ is a well-defined isometry is completely the same as in the case of von Neumann algebraic quantum groups, by the simple observation that $\Delta_N(x)^*\Delta_N(y) = \Delta_M(x^*y)$ for $x,y\in N$, and the fact that $(\psi_M\otimes \iota)(\Delta_M(x^*y)) = \psi_M(x^*y)1_M$ for $x,y\in \mathscr{N}_{N,\psi_M}$, by (polarization and the) definition of right-invariance.\\

\noindent But in this case, also the proof that $\widetilde{V}$ is a \emph{unitary} is easy. Indeed, since $\Delta_N(xy) = \Delta_N(x)\Delta_M(y)$ for $x\in N$ and $y\in M$, we have, for $x\in N$ and $y,z\in \mathscr{N}_{M,\psi_M}$, that \[\widetilde{V}\,\Gamma_N(xy)\otimes \Gamma_N(z) = \Delta_N(x)(\Gamma_M\otimes \Gamma_M)(\Delta_M(y)(1\otimes z)).\] Now elements of the form $(\Gamma_M\otimes \Gamma_M)(\Delta_M(y)(1\otimes z))$ have dense linear span in $\mathscr{L}^2(M)\otimes \mathscr{L}^2(M)$. Hence the range of $\widetilde{V}$ contains the closure of the set $\Delta_N(N)\cdot \mathscr{L}^2(M)\otimes \mathscr{L}^2(M)$. As $\Delta_N(N)(M\otimes M)$ is $\sigma$-weakly dense in $N\bar{\otimes} N$ by definition of a Galois co-object, we see that indeed the range of $\widetilde{V}$ equals $\mathscr{L}^2(N)\otimes \mathscr{L}^2(N)$, so that $\widetilde{V}$ is in fact a unitary.\\

\noindent Now we prove that $\widetilde{V} \in B(\mathscr{L}^2(N)) \bar{\otimes} N$. Using that $ B(\mathscr{L}^2(N))\bar{\otimes}N$ is a corner of $ B(\mathscr{L}^2(N))\bar{\otimes}Q$, it follows that it is sufficient to show that $(\omega\otimes\iota )(\widetilde{V}) \in N$ for each $\omega\in B(\mathscr{L}^2(N))_*$. We may further simplify by taking $\omega$ of the form $\langle\Gamma_{N}(z),\,\cdot\,\Gamma_N(y)\rangle$ for $y,z \in  \mathscr{N}_{N,\psi_M}$, as the linear span of such elements is dense in $B(\mathscr{L}^2(N))$. But then it follows from the definition of $\widetilde{V}$ and a Fubini type argument that \[ (\omega\otimes\iota )(\widetilde{V}) = (\psi_M\otimes \iota)((z^*\otimes 1)\Delta_N(y)) \in N,\] where we remark that $(z^*\otimes 1)\Delta_N(y)$ lies in the domain $\mathscr{M}_{(\iota\otimes \varphi_M)}$ of the operator valued weight $\iota\otimes \varphi_M$ from $Q\bar{\otimes} M$ to $Q=Q\otimes 1$, since $\Delta_N(y)^*\Delta_N(y)=\Delta_M(y^*y)$ and $(z^*z\otimes 1)$ are inside $\mathscr{M}_{(\iota\otimes \varphi_M)}^+$.\\

\noindent Finally, if $x \in  \mathscr{N}_{N,\psi_M}$ and $\omega\in N_* \subseteq Q_*$, we have the Cauchy-Schwarz inequality \[(\iota\otimes \omega)(\Delta_N(x))^*(\iota\otimes \omega)(\Delta_N(x)) \leq \|\omega\|\,(\iota\otimes |\omega|)(\Delta_M(x^*x)),\] where $|\omega|$ is the absolute value of $\omega$. It follows that $(\iota\otimes \omega)(\Delta_N(x)) \in  \mathscr{N}_{N,\psi_M}$. If there further exist $y \in  \mathscr{N}_{N,\psi_M}$ and $z\in \mathscr{N}_{M,\psi_M}$ such that $\omega$ is of the form $\langle\Gamma_{M}(z),\,\cdot\,\Gamma_N(y)\rangle$, it follows from the definition of $\widetilde{V}$ that \[(\iota\otimes \omega)(\widetilde{V})\Gamma_N(x) = \Gamma_N((\iota\otimes \omega)\Delta_N(x)).\] By the closedness of $\Gamma_N$ and the density of the linear span of such functionals in $N_*$, it follows that this formula holds for any $\omega\in N_*$.\\

\end{proof}

\begin{Def} Let $(N,\Delta_N)$ be a Galois co-object for a von Neumann algebraic quantum group $(M,\Delta_M)$. We call the unitary $\widetilde{V}$ the \emph{right regular $(N,\Delta_N)$-corepresentation} of $(N,\Delta_N)$. \end{Def}

\noindent Similarly, one can define a left such corepresentation $\widetilde{W}$, such that $\widetilde{W}^*$ will then be an element of $N\bar{\otimes} B(\mathscr{L}^2(N))$.\\

\noindent The following Proposition is an easy consequence of the definition of $\widetilde{V}$.

\begin{Prop} Let $(N,\Delta_N)$ be a Galois co-object for a von Neumann algebraic quantum group $(M,\Delta_M)$. Let $V$ be the regular right corepresentation for $(M,\Delta_M)$, and let $\widetilde{V}$ be the right regular $(N,\Delta_N)$-corepresentation for $(N,\Delta_N)$.
\begin{enumerate} \item For any $x\in N$, we have \[\widetilde{V}(x\otimes 1)\widetilde{V}^* = \Delta_N(x).\]
\item The following \emph{pentagonal equation} hold: \[\widetilde{V}_{12}\widetilde{V}_{13}V_{23} = \widetilde{V}_{23}\widetilde{V}_{12}.\]
\end{enumerate}
\end{Prop}

\begin{proof} Choose $y\in \mathscr{N}_{M,\psi_M}$ and $x\in N$. Then $xy\in \mathscr{N}_{N,\psi_N}$, and $\Gamma_N(xy)=x\Gamma_{M}(y)$. From this, it is immediately seen, using the definition of $\widetilde{V}$ and $V$, that $\widetilde{V}(x\otimes 1) = \Delta_N(x)V$, and hence $\widetilde{V}^*(x\otimes 1)V = \Delta_N(x)$.\\

\noindent Since $\widetilde{V}\in B(\mathscr{L}^2(N))\bar{\otimes}N$, and since we can implement $\Delta_N$ by $V$ and $\widetilde{V}$ by means of the first point, the pentagon identity for $\widetilde{V}$ can be rewritten as $(\iota\otimes\Delta_N)(\widetilde{V}) = \widetilde{V}_{12}\widetilde{V}_{13}$. It is then enough to prove that, for any $\omega_1,\omega_2\in N_*$, we have \[(\iota\otimes((\omega_1\otimes \omega_2)\circ \Delta_N))(\widetilde{V}) = (\iota\otimes \omega_1)(\widetilde{V})(\iota\otimes \omega_2)(\widetilde{V}).\] But this follows immediately by applying these operators to a vector $\Gamma_N(x)$ with $x\in \mathscr{N}_{N,\psi_M}$, and using the final part of the previous Proposition together with the coassociativity of $\Delta_N$.\\
\end{proof}

\section{Projective corepresentations of von Neumann bialgebras}

\noindent In order to prove Theorem \ref{TheovNa}, we will use the notion of a \emph{projective corepresentation} of a von Neumann bialgebra. This is not the most natural way of proving the Theorem, but the more direct manner would require a lot of the arguments which are very similar to the ones of \cite{DeC1}, some of which are quite technical and subtle. We therefore thought it better to avoid this, and to actually \emph{use} the results of \cite{DeC1}.\\

\noindent The notion of a projective corepresentation was already introduced in Definition \ref{DefProj}. Let us however state clearly here what we mean by an isomorphism between projective corepresentations.\\

\begin{Def} Let $(M,\Delta_M)$ be a von Neumann bialgebra. We call two projective corepresentations $\alpha_1$ and $\alpha_2$ of $(M,\Delta_M)$ on respective Hilbert spaces $\mathscr{H}_1$ and $\mathscr{H}_2$ \emph{unitary equivalent} if there exists an isomorphism $\gamma:B(\mathscr{H}_1)\rightarrow B(\mathscr{H}_2)$ such that $\alpha_2 = (\gamma\otimes \iota)\alpha_1$. \end{Def}

\noindent The crucial property of a projective corepresentation will be that it can be \emph{implemented}, in the same way as ordinary projective representations of a locally compact group can be implemented by choosing a (measurable) section $\mathcal{U}(H)/S^1\rightarrow \mathcal{U}(\mathscr{H})$, with $\mathcal{U}$ the (Polish) group of unitaries of a (separable) Hilbert space. The notion we need for this is the following.

\begin{Def} Let $(M,\Delta_M)$ be a von Neumann bialgebra, and $(N,\Delta_N)$ a (right) Galois co-object for $(M,\Delta_M)$. A \emph{(unitary) projective (left) $(N,\Delta_N)$-corepresentation} of $(M,\Delta_M)$ consists of a unitary $\mathcal{G}\in N\bar{\otimes} B(\mathscr{H})$ (i.e., unitary as a map from $\mathscr{L}^2(M)\otimes \mathscr{H}$ to $\mathscr{L}^2(N)\otimes \mathscr{H}$), satisfying the corepresentation property \[(\Delta_{N}\otimes \iota)\mathcal{G} = \mathcal{G}_{13}\mathcal{G}_{23}.\]

\noindent If $\mathcal{G}_1$ and $\mathcal{G}_2$ are two $(N,\Delta_N)$-corepresentations on respective Hilbert space $\mathscr{H}_1$ and $\mathscr{H}_2$, we call $\mathcal{G}_1$ and $\mathcal{G}_2$ unitary equivalent if there exists a unitary $u:\mathscr{H}_1\rightarrow \mathscr{H}_2$ such that $\mathcal{G}_2(1\otimes u) = (1\otimes u)\mathcal{G}_1$.\\

\noindent If $\mathscr{H}$ is a Hilbert space, $\alpha: B(\mathscr{H})\rightarrow M\bar{\otimes} B(\mathscr{H})$ a projective representation of $(M,\Delta_M)$ on $\mathscr{H}$, and $(N,\Delta_N)$ a Galois co-object for $(M,\Delta_M)$, we say that a projective $(N,\Delta_N)$-corepresentation $\mathcal{G}$ \emph{implements} $\alpha$ if \[\alpha(x) = \mathcal{G}^*(1\otimes x)\mathcal{G} \qquad \textrm{for all }x\in B(\mathscr{H}).\]
\end{Def}

\noindent It is easy to see that any projective $(N,\Delta_N)$-corepresentation $\mathcal{G}$ on a Hilbert space $\mathscr{H}$ implements in a unique way a projective corepresentation $\alpha$ on $\mathscr{H}$, precisely by the formula $\alpha(x) = \mathcal{G}^*(1\otimes x)\mathcal{G}$. The fact that this is a coaction follows immediately by the relation between $\mathcal{G}$ with $\Delta_N$.\\

 \noindent We next want to show that any projective corepresentation is implemented by an $(N,\Delta_N)$-projective corepresentation (for \emph{some} $(N,\Delta_N)$), but we first establish a uniqueness result. It will make use of the following Lemma.

\begin{Lem}\label{LemNoAnn} Let $(N,\Delta_N)$ be a Galois co-object for a von Neumann bialgebra $(M,\Delta_M)$, and let $\mathcal{G}$ be an $(N,\Delta_N)$-projective corepresentation on a Hilbert space $\mathscr{H}$. Then the $M$-linear span of the space $\{(\iota\otimes \omega)\mathcal{G}\mid \omega\in B(\mathscr{H})_*\}$ is $\sigma$-weakly dense in $N$.

\end{Lem}

\begin{proof} Let $\widetilde{N}$ be the $\sigma$-weak closure of the space $\{(\iota\otimes \omega)(\mathcal{G})m\mid \omega\in B(\mathscr{H})_*,m\in M\}$, and suppose that $\widetilde{N}\neq N$. Then, if $Q$ is the linking von Neumann algebra associated to $N$, there exists a non-zero annihilator $x\in Q_{11}$ of $\widetilde{N}$, again since this space is a non-trivial right $M$-submodule of $N$. But this means that $x(\iota\otimes \omega)(\mathcal{G})=0$ for all $\omega\in B(\mathscr{H})_*$, and hence $(x\otimes 1)\mathcal{G}=0$. As $\mathcal{G}$ is a unitary, we get $x=0$, a contradiction. Hence $\widetilde{N}= N$.

\end{proof}

\begin{Prop}\label{PropEqual} Let $(M,\Delta_M)$ be a von Neumann bialgebra, and $\alpha: B(\mathscr{H})\rightarrow M\bar{\otimes} B(\mathscr{H})$ a projective corepresentation of $(M,\Delta_M)$ on a Hilbert space $\mathscr{H}$.\\

\noindent Then if $(N_1,\Delta_{N_1})$ and $(N_2,\Delta_{N_2})$ are two Galois co-objects for $(M,\Delta_M)$, both equipped with a $(N_i,\Delta_{N_i})$-projective corepresentation $\mathcal{G}_i$ implementing $\alpha$, then there exists an isomorphism $\pi:(N_1,\Delta_{N_1})\rightarrow (N_2,\Delta_{N_2})$ of Galois co-objects such that $(\pi\otimes \iota)\mathcal{G}_1 = \mathcal{G}_2$.
\end{Prop}

\noindent Here, the notion of isomorphism for Galois co-objects is of course an isomorphism of Morita Hilbert W$^*$-modules intertwining the comultiplication structures.

\begin{proof} For $\xi,\eta$ vectors in $\mathscr{H}$, denote $\omega_{\xi,\eta}= \langle\xi,\,\cdot\,\eta\rangle$, and denote  $\theta_{\xi,\eta}$ for the finite rank operator $\chi \rightarrow \langle \eta,\chi\rangle \xi$. Then if $\xi_1,\xi_2,\eta_1,\eta_2$ are vectors in $\mathscr{H}$, it is easily seen that \begin{eqnarray*}(\iota \otimes \omega_{\xi_1,\eta_1})(\mathcal{G}_i)^*(\iota\otimes \omega_{\xi_2,\eta_2})(\mathcal{G}_i) &=& (\iota\otimes \omega_{\eta_1,\eta_2})(\mathcal{G}_i^*(1\otimes \theta_{\xi_1,\xi_2})\mathcal{G}_i) \\ &=& (\iota\otimes \omega_{\eta_1,\eta_2})(\alpha(\theta_{\xi_1,\xi_2})),\end{eqnarray*} for both $i\in \{1,2\}$.\\

\noindent The Proposition then follows immediately by the previous Lemma and Lemma \ref{LemuniMod2}.
\end{proof}

\noindent However, this does \emph{not} imply that a if $(N,\Delta_N)$ is a Galois co-object, and $\mathcal{G}_1$ and $\mathcal{G}_2$ two projective $(N,\Delta_N)$-corepresentations implementing the same projective corepresentation, that they are isomorphic. The reason is that for projective $(N,\Delta_N)$-corepresentations with \emph{fixed} $(N,\Delta_N)$, the notion of isomorphism is stronger. The concrete situation is the following.

\begin{Prop}\label{PropIsos} Let $(N,\Delta_N)$ be a Galois co-object for a von Neumann bialgebra, and let $(Q,\Delta_Q)$ be the associated linking weak von Neumann bialgebra. Suppose that $\mathcal{G}_1$ and $\mathcal{G}_2$ are two projective $(N,\Delta_N)$-corepresentations on a Hilbert space $\mathscr{H}$, such that \[\alpha(x) = \mathcal{G}_1^*(1\otimes x)\mathcal{G}_1 = \mathcal{G}_2^*(1\otimes x)\mathcal{G}_2\qquad \textrm{for all }x\in B(\mathscr{H}).\] Then there exists a group-like unitary $u\in P$ such that $\mathcal{G}_1 = (v\otimes 1)\mathcal{G}_2$.
\end{Prop}

\noindent We recall that the group-like property means that $\Delta_P(v) = v\otimes v$.

\begin{proof} As \[\mathcal{G}_1\mathcal{G}_2^* \in P\bar{\otimes} B(\mathscr{H})\subseteq B(\mathscr{L}^2(N))\bar{\otimes} B(\mathscr{H})\] commutes with all $(1\otimes x)$ with $x\in B(\mathscr{H})$, there exists a unitary $v\in P$ such that $\mathcal{G}_1 = (v\otimes 1)\mathcal{G}_2$. We then must show that $v$ is group-like. This follows by plugging in the above equality in the identities $(\Delta_N\otimes \iota)(\mathcal{G}_i) = (\mathcal{G}_i)_{13}(\mathcal{G}_i)_{23}$, using that $\Delta_N(xy)=\Delta_P(x)\Delta_N(y)$ for $x\in P$ and $y\in N$.
\end{proof}

\noindent Of course, it is still possible that $(v\otimes 1)\mathcal{G}$ and $\mathcal{G}$ are isomorphic, but this will not always be the case.\\

\noindent Let us now prove that any projective corepresentation is implemented. In \cite{DeC1}, we proved this for von Neumann algebraic quantum groups, but in a very roundabout way. Here, we will give a rather elementary proof which is valid in the more general setting of von Neumann bialgebras. Nevertheless, we will later on actually need the result as it appears in \cite{DeC1}, because it contains some more information.\\

\begin{Prop} Let $(M,\Delta_M)$ be a von Neumann bialgebra, $\mathscr{H}$ a Hilbert space, and $\alpha:B(\mathscr{H})\rightarrow M\bar{\otimes}B(\mathscr{H})$ a projective corepresentation of $(M,\Delta_M)$ on $\mathscr{H}$. Then there exists a Galois co-object $(N,\Delta_N)$ for $(M,\Delta_M)$, together with a projective $(N,\Delta_N)$-corepresentation $\mathcal{G}$ on $\mathscr{H}$ which implements $\alpha$.\end{Prop}

\begin{proof} Choose an index set $I$ with cardinality $\textrm{dim}(\mathscr{H})$, and let $0$ be a distinguished element of $I$. Choose a basis $\{e_i\mid i\in I\}$ of $\mathscr{H}$, and denote by $e_{ij}$ the matrix units in $B(\mathscr{H})$ with respect to this basis. Let further $\mathscr{K}$ be a Hilbert space on which $M$ is faithfully and normally represented, and denote $\mathscr{I} = \alpha(e_{00})(\mathscr{K}\otimes \mathscr{H})$.\\

\noindent We can then define a unitary \[\mathcal{G}: \mathscr{K}\otimes \mathscr{H} \rightarrow \mathscr{I} \otimes \mathscr{H} :\xi \rightarrow \sum_{i\in I}(\alpha(e_{0i})\xi)\otimes e_i,\] the adjoint being \[ \mathcal{G}^*: \mathscr{I}\otimes \mathscr{H}\rightarrow \mathscr{K}\otimes \mathscr{H}: \xi\otimes \delta_i\rightarrow \alpha(e_{i0})\xi.\] For any $x\in B(\mathscr{H})$, we have \[\mathcal{G}^*(1\otimes x)\mathcal{G} = \alpha(x),\] which follows most easily if one takes $x$ a matrix unit for example.\\

\noindent Denote then by $N$ the $\sigma$-weakly closed linear span of \[\{(\iota\otimes \omega_{0j})(\mathcal{G})m\mid j\in I,m\in M\} \subseteq B(\mathscr{K},\mathscr{I}).\] By definition, it is a right $M$-module. Moreover, just as in Proposition \ref{PropEqual} one has, denoting $\omega_{ij}= \langle e_i,\,\cdot\,e_j\rangle$, that \[(\iota \otimes \omega_{0j})(\mathcal{G})^*(\iota\otimes \omega_{0l})(\mathcal{G}) = (\iota\otimes \omega_{jl})(\alpha(e_{00})),\qquad i,j,k,l\in I.\] Hence $N$ becomes a Hilbert W$^*$-module by the formula $\langle x,y\rangle_M = x^*y$.\\

\noindent Now for $j,k,l\in I$, we have \begin{eqnarray*} (\iota\otimes \omega_{jk})(\mathcal{G}) &=& (\iota\otimes \omega_{0l})((1\otimes e_{0j})\mathcal{G}(1\otimes e_{kl}))\\ &=& (\iota\otimes \omega_{0l})(\mathcal{G}\alpha(e_{0j})(1\otimes e_{kl})) \\ &=& \sum_{i\in I} (\iota\otimes \omega_{0i})(\mathcal{G})(\iota\otimes \omega_{il})(\alpha(e_{0j})(1\otimes e_{kl})).\end{eqnarray*} Hence, each $(\iota\otimes \omega_{jk})(\mathcal{G})$ lies in $N$, and thus $\mathcal{G}\in N\bar{\otimes}B(\mathscr{H})$.\\

\noindent In particular, we have $(\iota\otimes \omega_{i0})(\mathcal{G})\in N$. As \[(\iota \otimes \omega_{i0})(\mathcal{G})^*(\iota\otimes \omega_{i0})(\mathcal{G}) = (\iota\otimes \omega_{00})(\alpha(e_{ii}))\] for $i\in I$, we see that the linear span of the range of $\langle \,\cdot\,,\,\cdot\,\rangle_M$ contains $(\iota\otimes \omega_{00})(\alpha(1)) = 1_M$, and so $N$ is a full Hilbert $M$-module. As it arises as a $\sigma$-weakly closed subspace of $B(\mathscr{K},\mathscr{I})$, we have that $N$ is a Morita Hilbert $M$-module by Proposition \ref{PropTern}.\\

\noindent Denote \[x_{i} = (\iota\otimes \omega_{0i})(\mathcal{G})\in N,\] and \[y_{i} = (\iota\otimes \iota\otimes\omega_{0i})(\mathcal{G}_{13}\mathcal{G}_{23}) \in N\bar{\otimes} N,\] with $i\in I$. As $x_{i}^*x_{j} = (\iota\otimes \omega_{ij})(\alpha(e_{00}))$, while \begin{eqnarray*} y_{i}^*y_{j} &=& (\iota\otimes \iota\otimes \omega_{ij})((\iota\otimes \alpha)\alpha(e_{00})) \\ &=& \Delta_M(x_{i}^*x_{j})\end{eqnarray*} by an easy computation, we can apply Lemma \ref{LemuniMod2} to obtain a $\Delta_M$-compatible morphism $\Delta_N: N\rightarrow N\bar{\otimes}N$ such that $(\Delta_N\otimes \iota)\mathcal{G} = \mathcal{G}_{13}\mathcal{G}_{23}$. The Proposition will then be proven if we can show that $(N,\Delta_N)$ is a Galois co-object.\\

\noindent In fact, by the above compatibility with $\mathcal{G}$, and the fact that the first leg of $\mathcal{G}$ generates $N$ as a right $M$-module, it follows immediately that $\Delta_N$ will be coassociative. The only thing which remains then is to see if $\Delta_N(N)(M\bar{\otimes}M)$ is $\sigma$-weakly dense in $N\bar{\otimes}N$. But this follows precisely as in the proof of Lemma \ref{LemNoAnn}.
\end{proof}

\noindent \emph{Remark:} In particular, the foregoing allows one to construct from a projective corepresentation of $(M,\Delta_M)$ (i.e.~ a coaction on a type $I$-factor) a Galois co-object $(N,\Delta_N)$, and hence, by Proposition \ref{PropCooblink}, a linking weak von Neumann bialgebra $(Q,\Delta_Q)$, which contains in turn a (possibly) \emph{new} von Neumann bialgebra $(P,\Delta_P)$ in its upper left corner. In \cite{DeC5}, we applied this construction to the action of $SU_q(2)$ on the standard Podle\`{s} sphere (whose associated von Neumann algebra is indeed a type $I$-factor) to `rediscover' Woronowicz' quantum $E(2)$ group (\cite{Wor2}). In \cite{DeC6}, we applied it to the action of $SU_q(2)$ on a $\mathbb{Z}_2$-quotient of the equatorial Podle\'{s} sphere (which can be interpreted as a quantized projective plane, with again a type $I$-factor as its associated von Neumann algebra) to `rediscover' the extended quantum $SU(1,1)$ group (as it appears in \cite{Koe1}). We hope in future work to obtain in this way some interesting $q$-deformations of higher-dimensional non-compact Lie groups.\\

\noindent The following Proposition will be an immediate corollary of Proposition \ref{PropEqual} and the results of \cite{DeC1}. We first remark however that the object $(\widehat{Q},\Delta_{\widehat{Q}})$ which appears in the beginning of the first section of \cite{DeC1} is a linking weak von Neumann bialgebra in the sense of the present paper. Indeed, $\widehat{Q}$ is a linking von Neumann algebra by the remark following Proposition \ref{PropLinkk} of the present paper, and since $\Delta_{\widehat{Q}}$ was \emph{constructed} in \cite{DeC1} as a unital map $\widehat{Q}\rightarrow \widehat{Q}* \widehat{Q}$, it will hence make $\widehat{Q}$ a linking weak von Neumann bialgebra. In fact, as we showed in \cite{DeC1} that this specific linking von Neumann bialgebra has von Neumann algebraic quantum groups at its corners, it is a von Neumann algebraic linking quantum groupoid in the terminology of the present paper.

\begin{Prop}\label{PropvNa} Let $\alpha$ be a projective corepresentation of a von Neumann algebraic quantum group $(M,\Delta_M)$ on a Hilbert space $\mathscr{H}$. Let $(N,\Delta_N)$ be a Galois co-object for which there exists a projective $(N,\alpha)$-corepresentation implementing $\alpha$. Then the linking weak von Neumann bialgebra associated to $(N,\Delta_N)$ is a von Neumann algebraic linking quantum groupoid.
\end{Prop}

\begin{proof} In Theorem 6.2 of \cite{DeC1}, we constructed a von Neumann algebraic linking quantum groupoid $(\widetilde{Q},f,\Delta_{\widetilde{Q}})$, such that the Galois co-object $(\widetilde{Q}_{12},\Delta_{12})$ had a projective $(\widetilde{Q}_{12},\Delta_{12})$-corepresentation implementing $\alpha$. If then $(Q,e,\Delta_Q)$ is the linking weak von Neumann bialgebra associated to $(N,\Delta_N)$, we have, by Proposition \ref{PropEqual} and Lemma \ref{Lemid}, an isomorphism $\pi$ from $(Q,e,\Delta_Q)$ to $(\widetilde{Q},f,\Delta_{\widetilde{Q}})$, which intertwines $\Delta$ because the restriction to the 12-part does. Hence $(Q_{11},\Delta_{11})$ is a von Neumann algebraic quantum group.
\end{proof}

\noindent Finally, we use the previous Proposition to prove Theorem \ref{TheovNa}.\\

\begin{proof}[Proof (of Theorem \ref{TheovNa}).] Let $(M,\Delta_M)$ be a von Neumann algebraic quantum group, and let $(N,\Delta_N)$ be a Galois co-object for $(M,\Delta_M)$. By Proposition \ref{PropvNa}, it is enough to show that there exists a left $(N,\Delta_N)$-corepresentation. But it is easy to see that, with $\Sigma$ denoting the flip map and $\widetilde{V}$ the regular right $(N,\Delta_N)$-corepresentation associated with $(N,\Delta_N)$, we have that $\Sigma\widetilde{V}\Sigma \in N\bar{\otimes}B(\mathscr{L}^2(N))$ is a left $(N,\Delta_N)$-corepresentation by the pentagonal equation for $\widetilde{V}$. This concludes the proof.
\end{proof}

\noindent \emph{Remark:} The above proof is of course very sparse with information on how the invariant weights on the comonoidally W$^*$-Morita equivalent von Neumann bialgebra $(P,\Delta_P)$ are obtained. The crucial point to observe is that in \cite{DeC1}, we proved that there exists a one-parameter-group of unitaries on $\mathscr{L}^2(N)$ which implements the modular one-parametergroup (of say the left invariant weight) on $\pi_r(M)$, where $\pi_r$ is the natural right representation of $M$ on $\mathscr{L}^2(N)$. A theorem due to Connes implies that this one-parameter-family is in fact generated by the spatial derivative between (the opposite of) the left invariant weight on $M$ and a \emph{uniquely determined} weight on $P$. We then showed that this new weight is left invariant.\\

\noindent The way in which the above-mentioned one-parameter-group of unitaries was constructed is in itself not so straightforward, and is heavily influenced by the way in which all structures on a von Neumann algebraic quantum group interact with each other. In any case, even though the intuition from \cite{DeC1} could in principle be used to prove Theorem \ref{TheovNa} without recourse to the (dual) theory in \cite{DeC1}, we have deemed this task not worth the effort, as there seemed to be little gain in reiterating all technical arguments.

\section{2-cocycles}

\noindent Let us now briefly consider the special case of \emph{cleft} Galois co-objects, which are those Galois co-objects constructed from a unitary \emph{2-cocycle} (\cite{Eno1}). This discussion will then supplement the one in the fifth section of \cite{DeC1}.

\begin{Def} Let $(M,\Delta_M)$ be a von Neumann bialgebra, and $\Omega\in M\bar{\otimes}M$ a unitary. We call $\Omega$ a \emph{unitary 2-cocycle} if $\Omega$ satisfies the \emph{2-cocycle identity}: \[(\Omega\otimes 1)(\Delta_M\otimes \iota)(\Omega) = (1\otimes \Omega)(\iota\otimes \Delta_M)(\Omega).\] \end{Def}

\begin{Prop} If $\Omega$ is a unitary 2-cocycle for a von Neumann algebraic quantum group $(M,\Delta_M)$, then \[(N,\Delta_N) := (M,\Omega\Delta_M(\,\cdot\,))\] with $\langle x,y\rangle_M= x^*y$ for $x,y\in M$, is a right Galois co-object for $(M,\Delta_M)$.
\end{Prop}

\begin{proof} The fact that $\Delta_N$ is coassociative is immediate from the 2-cocycle identity. Also the other properties of Galois co-objects are trivial to verify.\end{proof}

\noindent The following Propositions are quite trivial to prove, but it is important to note them.

\begin{Prop} If $(N,\Delta_N)$ is a Galois co-object for a von Neumann algebraic quantum group, and $N\cong M$ as right Hilbert W$^*$-modules, then there exists a unitary 2-cocycle $\Omega$ such that $(N,\Delta_N)\cong (M,\Omega\Delta_M(\,\cdot\,))$.\end{Prop}

\begin{proof} Identifying $N$ with $M$ as a right Hilbert W$^*$-module, we have that $\Omega = \Delta_N(1_M)$ is a unitary, satisfying the 2-cocycle condition since $\Delta_N$ is coassociative and $\Delta_N(x)=\Delta_N(1_M)\Delta_M(x)$ for $x\in M$. This final identity then also proves that $\Delta_N=\Omega\Delta_M(\,\cdot\,)$.

\end{proof}

\noindent Hence these Galois co-objects can be characterized as those for which the associated underlying W$^*$-Morita equivalence (i.e.~ without the comonoidal structure) is trivial.

\begin{Prop} Let $\Omega_1$ and $\Omega_2$ be two unitary 2-cocycles for a von Neumann algebraic quantum group $(M,\Delta_M)$, and let $(N_1,\Delta_{N_1})$ and $(N_2,\Delta_{N_2})$ be the associated Galois co-objects. Then $(N_1,\Delta_{N_1})$ and $(N_2,\Delta_{N_2})$ are isomorphic iff $\Omega_1$ and $\Omega_2$ are coboundary equivalent, in that there exists a unitary $v\in M$ such that \[\Omega_2 = (v^*\otimes v^*)\Omega_1\Delta_M(v).\]

\end{Prop}

\begin{proof} If $\Omega_1$ and $\Omega_2$ are coboundary equivalent by a unitary $v$, it is immediately verified that left multiplication by $v^*$ provides an isomorphism between $(N_1,\Delta_{N_1})$ and $(N_2,\Delta_{N_2})$.\\

\noindent Conversely, suppose that $(M,\Omega_1\Delta_M(\,\cdot\,))$ and $(M,\Omega_2\Delta_M(\,\cdot\,))$ are isomorphic as right $N$-Galois co-objects by a map $\phi$. Then $\phi(1_M)$ is a unitary, whose adjoint we denote by $v$. Then $\phi(m) = v^*m$ for all $m\in M$. As $\phi$ intertwines the coproducts, we find that $\Omega_2 \Delta_M(v^*)= (v^*\otimes v^*)\Omega_1$, so that $\Omega_1$ and $\Omega_2$ are coboundary equivalent.

\end{proof}

\end{document}